\documentclass[conference,a4paper]{IEEEtran}
\IEEEoverridecommandlockouts
\usepackage{cite}
\usepackage{amsmath}    
\usepackage{amssymb,amsthm}    
\usepackage{graphicx}
\usepackage{textcomp}
\usepackage{xcolor}

\usepackage{chemformula}
\usepackage{bm}
\usepackage{subcaption}
\usepackage{siunitx}
\usepackage{algorithm}
\usepackage{algpseudocode}
\usepackage{soul}
\usepackage{hyperref}

\newcommand{\T}{\mathsf{T}}
\renewcommand{\L}{\mathcal{L}}
\renewcommand{\b}{\boldsymbol}

\newcommand{\x}{\b{x}}

\graphicspath{{./figures/}}

\def\BibTeX{{\rm B\kern-.05em{\sc i\kern-.025em b}\kern-.08em
		T\kern-.1667em\lower.7ex\hbox{E}\kern-.125emX}}
\begin{document}

\title{Implicit-Explicit Error Indicator based on Approximation Order\
	\thanks{The authors would like to acknowledge the financial support of
		Slovenian Research Agency (ARRS) in the framework of the research core funding
		No. P2-0095.}
}

\makeatletter
\newcommand{\linebreakand}{%
	\end{@IEEEauthorhalign}
	\hfill\mbox{}\par
	\mbox{}\hfill\begin{@IEEEauthorhalign}
}
\makeatother

\author{
	\IEEEauthorblockN{Mitja Jan\v{c}i\v{c}}
	\IEEEauthorblockA{\textit{Parallel and Distributed Systems Laboratory} \\
		\textit{Jo\v{z}ef Stefan Institute}\\
		\textit{Jo\v{z}ef Stefan International Postgraduate School} \\
		Ljubljana, Slovenia \\
		mitja.jancic@ijs.si}
	\and
	\IEEEauthorblockN{Filip Strni\v{s}a}
	\IEEEauthorblockA{\textit{Parallel and Distributed Systems Laboratory} \\
		\textit{Jo\v{z}ef Stefan Institute}\\
		Ljubljana, Slovenia \\
		filip.strnisa@ijs.si}
	\linebreakand 
	\centering
	\IEEEauthorblockN{Gregor Kosec}
	\IEEEauthorblockA{\textit{Parallel and Distributed Systems Laboratory} \\
		\textit{Jo\v{z}ef Stefan Institute}\\
		Ljubljana, Slovenia \\
		gregor.kosec@ijs.si}
}

\maketitle

\begin{abstract}
	With the immense computing power at our disposal, the numerical solution of partial differential equations (PDEs) is becoming a day-to-day task for modern computational scientists. However, the complexity of real-life problems is such that tractable solutions do not exist. This makes it difficult to validate the numerically obtained solution, so good error estimation is crucial in such cases. It allows the user to identify problematic areas in the computational domain that may affect the stability and accuracy of the numerical method. Such areas can then be remedied by either \textit{h}- or \textit{p}-adaptive procedures.
	In this paper, we propose to estimate the error of the numerical solution by solving the same governing problem implicitly and explicitly, using a different approximation order in each case. We demonstrate the newly proposed error indicator on the solution of a synthetic two-dimensional Poisson problem with tractable solution for easier validation. We show that the proposed error indicator has good potential for locating areas of high error.
\end{abstract}

\bigskip


\section{Introduction}
In physical modelling, systems of partial differential equations (PDEs) are used to describe the dynamical properties of many natural phenomena. Moreover, the solution of such systems is often of interest to engineers and scientists. However, due to their complexity, they almost never have analytical solutions, and need to be treated numerically, leading to a numerical solution.
In general, PDE problems are often solved using one of the following three methods: the finite volume method (FVM), the finite element method (FEM) and the finite difference method (FDM).
Recently, however, a generalised formulation of FDM, the radial basis function-generated finite differences (RBF-FD)~\cite{tolstykh2000using, tolstykh2003using}, has become increasingly popular. This is mainly because RBF-FD is a variant of the mesh-free methods ~\cite{belytschko1996meshless}, i.e. the method can operate on scattered nodes, unlike the previously mentioned mesh-based methods.

In the context of RBF-FD, linear differential operators are approximated over a set of RBFs augmented with monomials. Augmentation is necessary to ensure convergent and stable behaviour of the method~\cite{bayona2017role,flyer2016}. Additionally, it also enables a direct control over the order of the approximation method, as it corresponds to the highest order used in the approximation basis.

Nevertheless, after the numerical solution is obtained, scientists are often confronted with the difficulty of validating it. For that reason, researchers proposed error indicators~\cite{9597066, carstensen2005explicit} to identify problematic areas with a high error of the numerical solution. In practise, different adaptive numerical methods are then applied to these areas~\cite{SEGETH20101589} ensuring a finer local field description (\textit{h}-adaptivity) or higher polynomial degree approximations (\textit{p}-adaptivity), effectively improving the accuracy of numerical solution.

In this paper, we present an \textit{a posteriori} error indicator that measures the error of an implicit solution. The error indicator is applied through the meshless RBF-FD method as found in the \emph{Medusa library}~\cite{slak2021}. In general, the idea is to apply higher order explicit differential operators approximations to the implicitly obtained solution and thus indicate the areas with high error of the numerical solution. In the continuation of this work, the proposed error indicator will be named IMEX (\emph{im}plicit-\emph{ex}plicit) error indicator.

\section{IMEX error indicator}
Let there be a partial differential equation of type:
\begin{equation}
	\label{eq: general PDE}
	\L u = a,
\end{equation}
where $\L$ is an arbitrary partial differential operator applied to $u$, and equaling the constant $a$.
Such a problem is first solved implicitly, using a lower-order approximation of $\L$, $\L^{(lo)}$, obtaining the solution $u^{(im)}$ in the process.
The $u^{(im)}$ is then used to reconstruct $a$ explicitly with the help of higher-order approximation of $\L$, $\L^{(hi)}$, giving $a^{(ex)}$.
Finally, $a^{(ex)}$ is then tested against the analytical $a$ to indicate the error.
These steps can be summarized as follows:
\begin{enumerate}
	\item compute approximations $\L^{(lo)}$ and $\L^{(hi)}$;
	\item solve $\L^{(lo)} u = a$ implicitly, obtain $u^{(im)}$;
	\item compute $a^{(ex)} = \L^{(hi)} u^{(im)}$;
	\item compare $a^{(ex)}$ and $a$ to indicate high error areas.
\end{enumerate}

\section{RBF-FD approximation of differential operators}
\label{sec:rbffd}

Since the introduction of meshless methods in the 1970s, many variants have been proposed. The first mention of RBF-FD dates from 2000 with the introduction from Tolstykh~\cite{tolstykh2000using}. Since then, the method has been thoroughly studied and applied to many real-world problems with recent applications to fluid flow~\cite{rot} and plasticity~\cite{filip} problems.

In the framework of RBF-FD, a linear differential operator $\L$ in the node $\x_c$ is approximated over a set of $n$ neighbouring (often called \emph{stencil}) nodes
\begin{eqnarray}
	\label{eq:approx}
	\widehat{\L u}(\x_c)=\sum_{i=1}^nw_iu(\x_i)
\end{eqnarray}
for an arbitrary function $u$ and weights $\b w$ yet to be determined. The weights $\b w$ are obtained by constructing a localised RBF approximation with a given set of radial basis functions (RBFs) $\theta$ centred at the stencil nodes of a central node $\x _c$
\begin{eqnarray}
	\theta(\x) = \theta(\left\| \x - \x_c\right \|).
\end{eqnarray}
The localized intepolation~\eqref{eq:approx} can be written in a linear system
\begin{equation}
	\underbrace{
		\begin{bmatrix}
			\theta(\x_1) & \cdots & \theta(\x_1) \\
			\vdots       & \ddots & \vdots       \\
			\theta(\x_n) & \cdots & \theta(\x_n) \\
		\end{bmatrix}
	}_{\b \Theta}
	\underbrace{
		\begin{bmatrix}
			w_1    \\
			\vdots \\
			w_n    \\
		\end{bmatrix}}_{\b w} =
	\underbrace{
		\begin{bmatrix}
			(\L \theta_1(\x)\big|_{ \b x = \b x_c} \\
			\vdots                                 \\
			(\L \theta_n(\x)\big|_{ \b x = \b x_c} \\
		\end{bmatrix}}_{\ell_\theta}.
\end{equation}

However, as previously observed by Bayona et al.~\cite{bayona2017role}, RBFs alone do not guarantee convergent behaviour or solvability of the system. To mitigate these problems, the approximation basis is extended by a set of $s = \binom{m+d}{d}$ monomials with up to and including degree $m$ in a $d$-dimensional domain.

With the additional constraints, the RBF-FD approximation can be written compactly as
\begin{equation}
	\label{eq:rbf-system-aug}
	\begin{bmatrix}
		{\b \Theta} & {\b P} \\
		{\b P}^\T   & \b 0
	\end{bmatrix}
	\begin{bmatrix}
		\b w \\
		\b \lambda
	\end{bmatrix}
	=
	\begin{bmatrix}
		\b\ell_\theta \\
		\b\ell_p
	\end{bmatrix},
\end{equation}
where $\b P$ is a $n\times s$ matrix of monomials evaluated at stencil points, $\b \ell_p$ is the vector of values composed by applying the operator under consideration $\L$ to the polynomials at $\x_c$, i.e.\ $\ell_p^i = (\L p_i(\x))\big|_{ \b x = \b x_c}$ and $\b \lambda $ are Lagrangian multipliers (which we discard after the solution had been obtained).

\section{Example}
The IMEX error indicator's performance is demonstrated on a synthetic example, which is commonly used when testing adaptive algorithms in mesh-based methods~\cite{mitchell2013}.

The example is the Poisson equation, which is solved in a 2D circular domain $\Omega$ with its center at (0, 0), and radius 1:
\begin{equation}
	\label{eq: poisson}
	\begin{aligned}
		\nabla^2 u                             & = f_{lap}(\bm{x}) \ \     & \mathrm{in} \ \Omega,                     \\
		\frac{\mathrm{d} u}{\mathrm{d} \bm{n}} & = \b f_{neu}(\bm{x}) \ \  & \mathrm{on} \ \partial\Omega, \ x \leq 0, \\
		u                                      & = f_{dir}(\bm{x}) \ \     & \mathrm{on} \ \partial\Omega, \ x > 0.
	\end{aligned}
\end{equation}
The Neumann, and Dirichlet boundary conditions are defined through $f_{neu}$, and $f_{dir}$, respectively:
\begin{equation}
	\label{eq: fneu}
	\b f_{neu}(\bm{x}) = -2 \alpha {\left[\exp{\left(-\alpha ||\bm{x} - \bm{x}_s||^2 \right)}\right]} \bm{x},
\end{equation}
\begin{equation}
	\label{eq: dir}
	f_{dir}(\bm{x}) = \exp{\left(-\alpha ||\bm{x} - \bm{x}_s||^2 \right)}.
\end{equation}
From these one can derive the analytical solution of the Laplacian $f_{lap}$ at point $\bm{x} = (x, y)$:
\begin{equation}
	\label{eq: flap}
	f_{lap}(\bm{x}) = 4 {\left(\alpha^2 ||\bm{x} - \bm{x}_s||^2 - \alpha\right)} \exp{\left(-\alpha ||\bm{x} - \bm{x}_s||^2 \right)}.
\end{equation}
The source is positioned at $\bm{x}_s$ while $\alpha$ controls the source strength.
$\bm{n}$ is the boundary normal at $\bm{x}$ on $\partial \Omega$.
$\bm{x}_s$ is positioned at (0.5, 0.5), and $\alpha$ is set to 1000.

The example was solved on a laptop with Intel Core i7-8750H CPU, and 16 GB RAM.
Results were computed, and written into a file in about 2 s\footnote{The source code for the example can be found at: \url{https://gitlab.com/e62Lab/2022_CP_splitech_IMEX_error_indicator_poisson_eg}}.

\section{Results and Discussion}
The computational domain is discretized and filled with scattered nodes using Medusa's built-in algorithms \cite{slak2019,slak2021}.
This procedure results in a domain discretized with 24882 points.
An example solution is shown in Fig.~\ref{fig: domain}.
\begin{figure}
	\centering
	\includegraphics[width=\columnwidth]{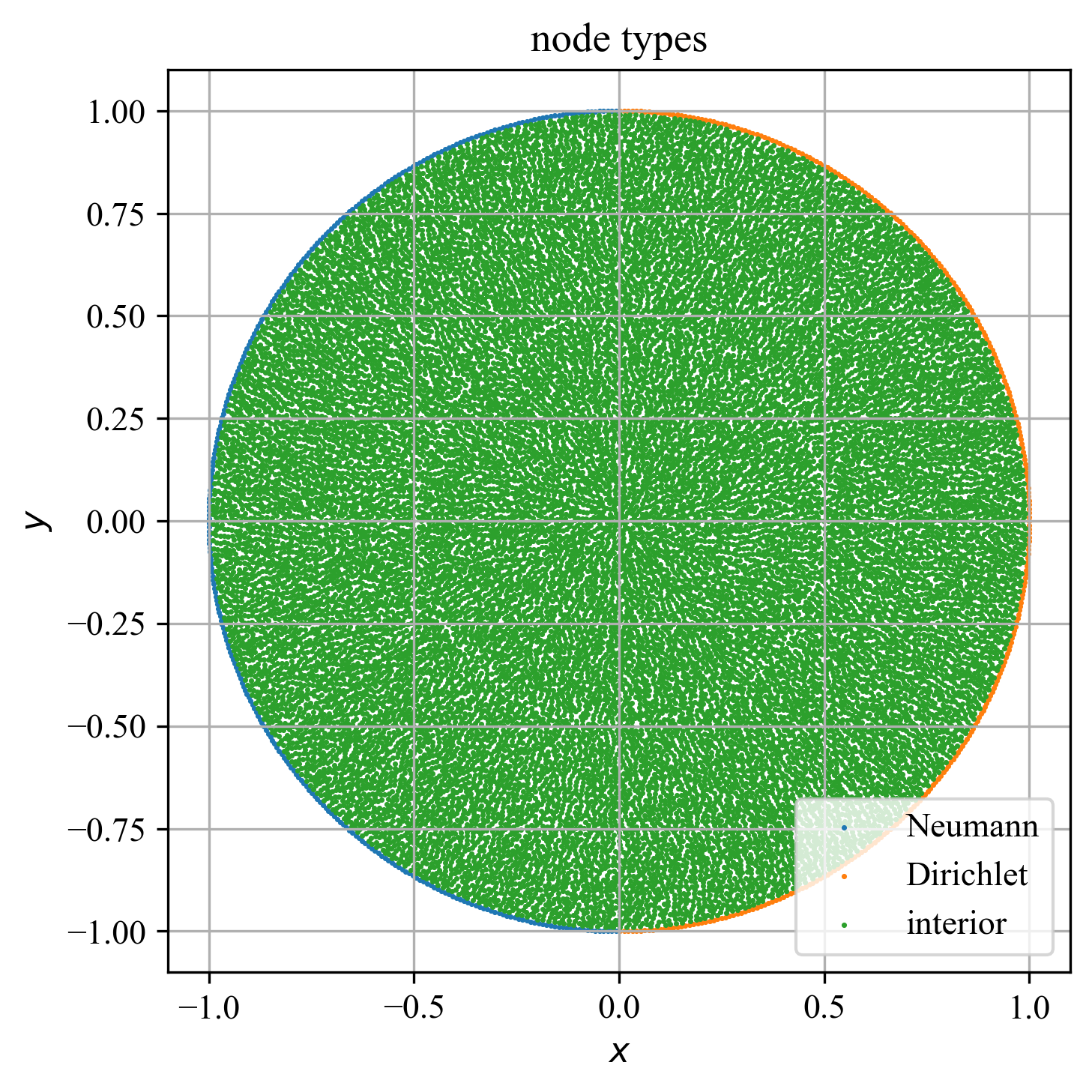}
	\caption{Domain discretization displaying positions of Neumann, and Dirichlet boundaries, as well as interior nodes.}\label{fig: domain}
\end{figure}
Support sizes for $\L^{(lo)}$, and $\L^{(hi)}$ are set to $2\binom{m+d}{d}$ (following the recommendations by Bayona et al.~\cite{bayona2017role}), $m$ being the monomial degree, and $d$ the number of dimensions of the domain.
The system in Eq.~\eqref{eq: poisson} is first solved implicitly, with lower order approximation of differential operators ${\nabla^2}^{(lo)}$, and $\frac{\mathrm{d}}{\mathrm{d} \bm{n}}^{(lo)}$, which were obtained with 2\textsuperscript{nd} degree monomials.
The solution for the scalar field $u^{(im)}$ is obtained with Eigen's \verb*|BiCGSTAB| solver \cite{eigenweb}.
To compute the RHS explicitly, a higher order approximation of the operator ${\nabla^2}^{(hi)}$, obtained with 4\textsuperscript{th} degree monomials, is applied to $u^{(im)}$.
The results are then compared to produce the IMEX error indicator $\epsilon_{IMEX}$:
\begin{equation}
	\label{eq: imex err}
	\epsilon_{IMEX} = \left|{\nabla^2}^{(hi)} u^{(im)}(\bm{x}) - f_{lap}(\bm{x})\right|.
\end{equation}
\begin{figure}
	\centering
	\includegraphics[width=\columnwidth]{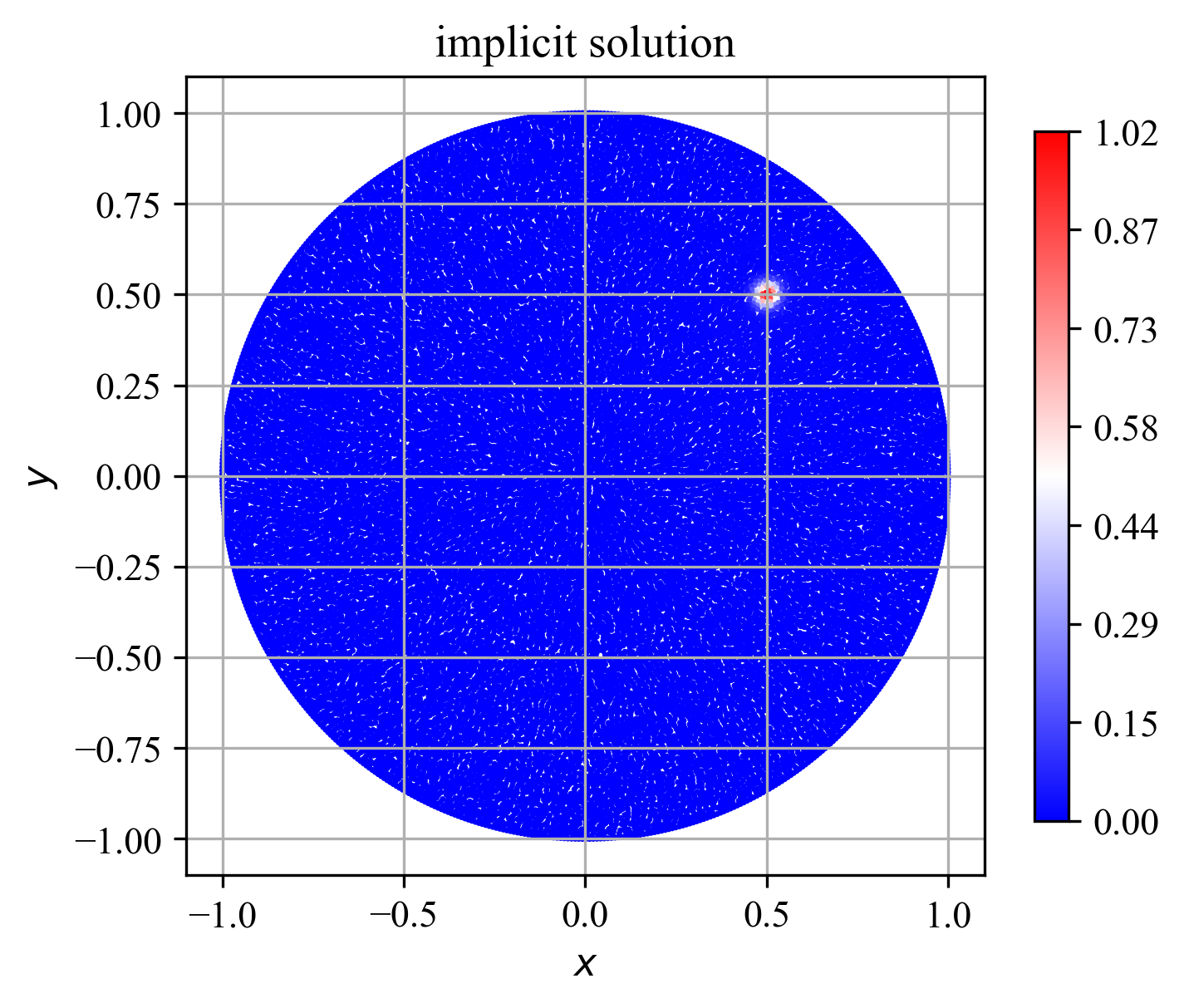}
	\caption{An example of the implicit solution.}\label{fig: solution}
\end{figure}
For validation purposes, the error of $u^{(im)}$, $\epsilon_{an}$, is also computed by comparing the implicit to the analytical solution. The latter is obtained with Eq.~\eqref{eq: dir}, and $\epsilon_{an}$ is:
\begin{equation}
	\label{eq: an err}
	\epsilon_{an} = \left|u^{(im)}(\bm{x}) - f_{dir}(\bm{x})\right|.
\end{equation}
\begin{figure}
	\centering
	\includegraphics[width=\columnwidth]{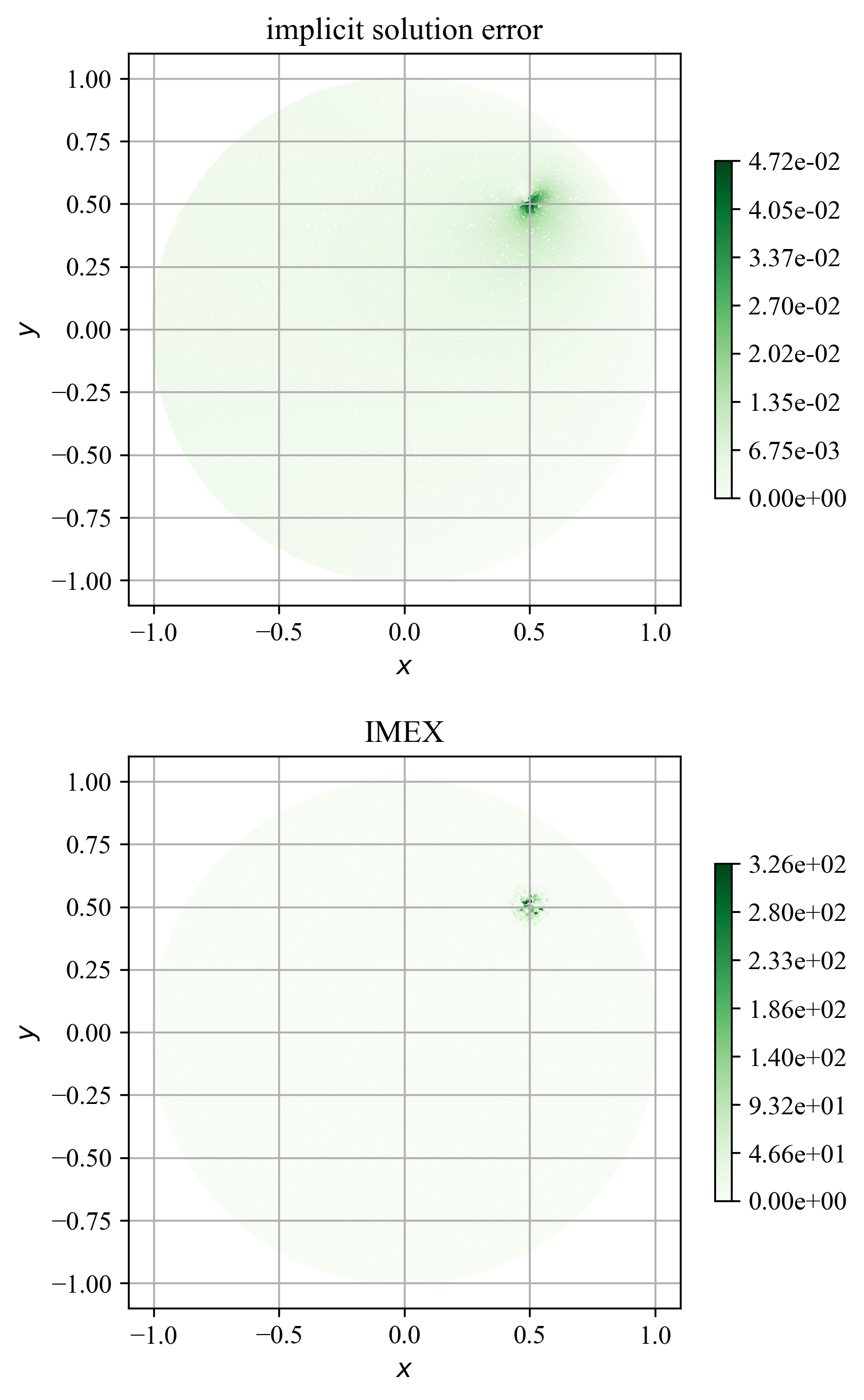}
	\caption{Comparison plots of $\epsilon_{an}$ (above), and $\epsilon_{IMEX}$ (below).}\label{fig: imex}
\end{figure}
\begin{figure}
	\centering
	\includegraphics[width=\columnwidth]{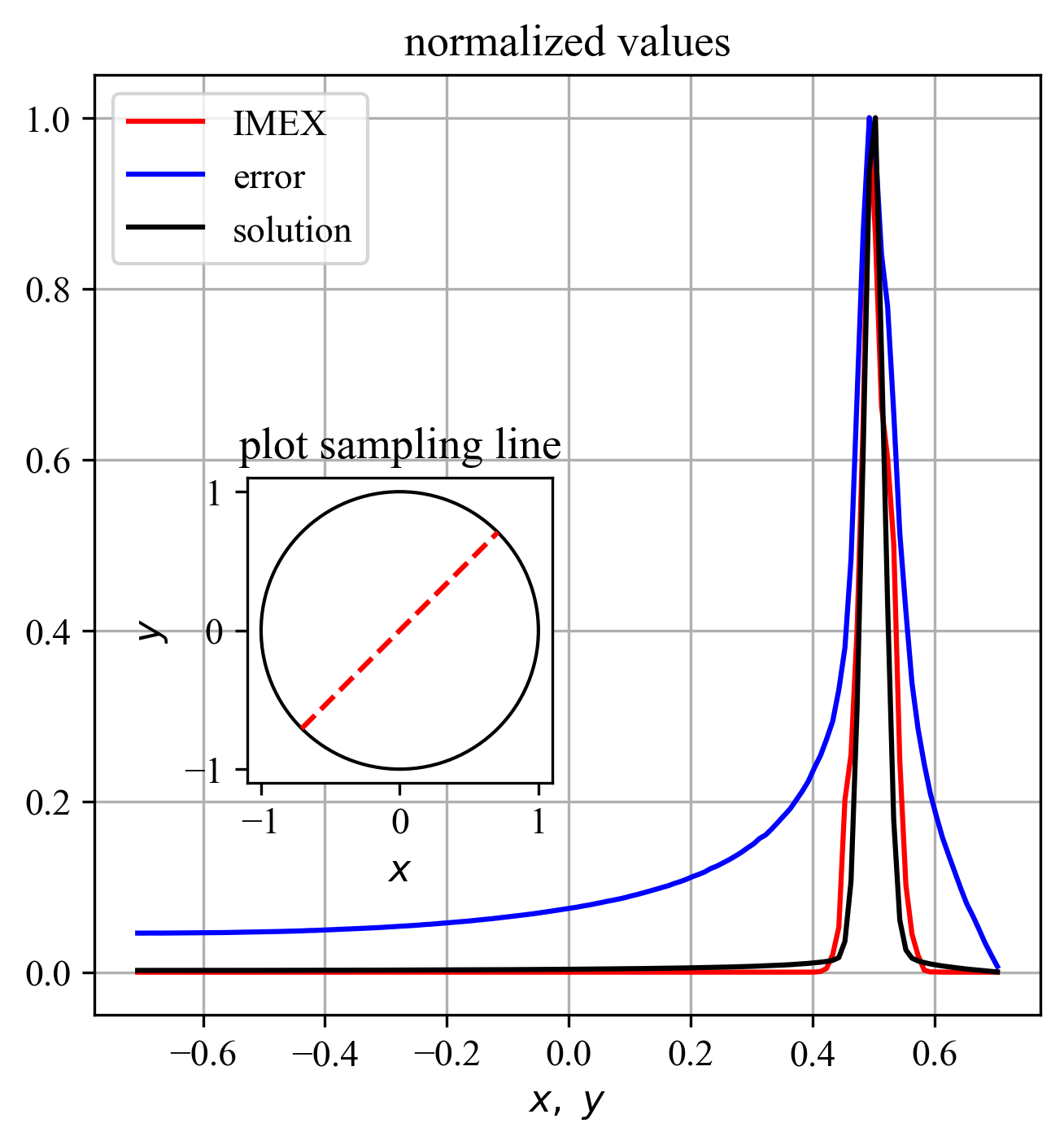}
	\caption{Implicit solution (solution), $\epsilon_{an}$ (error), and $\epsilon_{IMEX}$ (IMEX) normalized to their respective maximal values, plotted along the line $y = x; \ x,y \in \Omega$, $\L^{(hi)}$ are computed with 4\textsuperscript{th} degree monomials.}\label{fig: lines}
\end{figure}

Fig.~\ref{fig: solution} is displaying the implicit solution $u^{(im)}$ of Eq.~\eqref{eq: general PDE}, while $\epsilon_{an}$ and $\epsilon_{IMEX}$ are plotted in Fig.~\ref{fig: imex}.
For better clarity the implicit solution, $\epsilon_{an}$, and $\epsilon_{IMEX}$ are plotted in Fig.~\ref{fig: lines} along the line $y = x; \ x,y \in \Omega$.
As the solution was obtained on scattered nodes, the source for the aforementioned line is obtained by Shepard interpolation (Python, \verb*|ShepardIDWInterpolator| from \verb*|photutils.utils|~\cite{larry_bradley_2020_4044744}), sampling each plot line point from 9 nearest neighbors.
\begin{figure}
	\centering
	\includegraphics[width=\columnwidth]{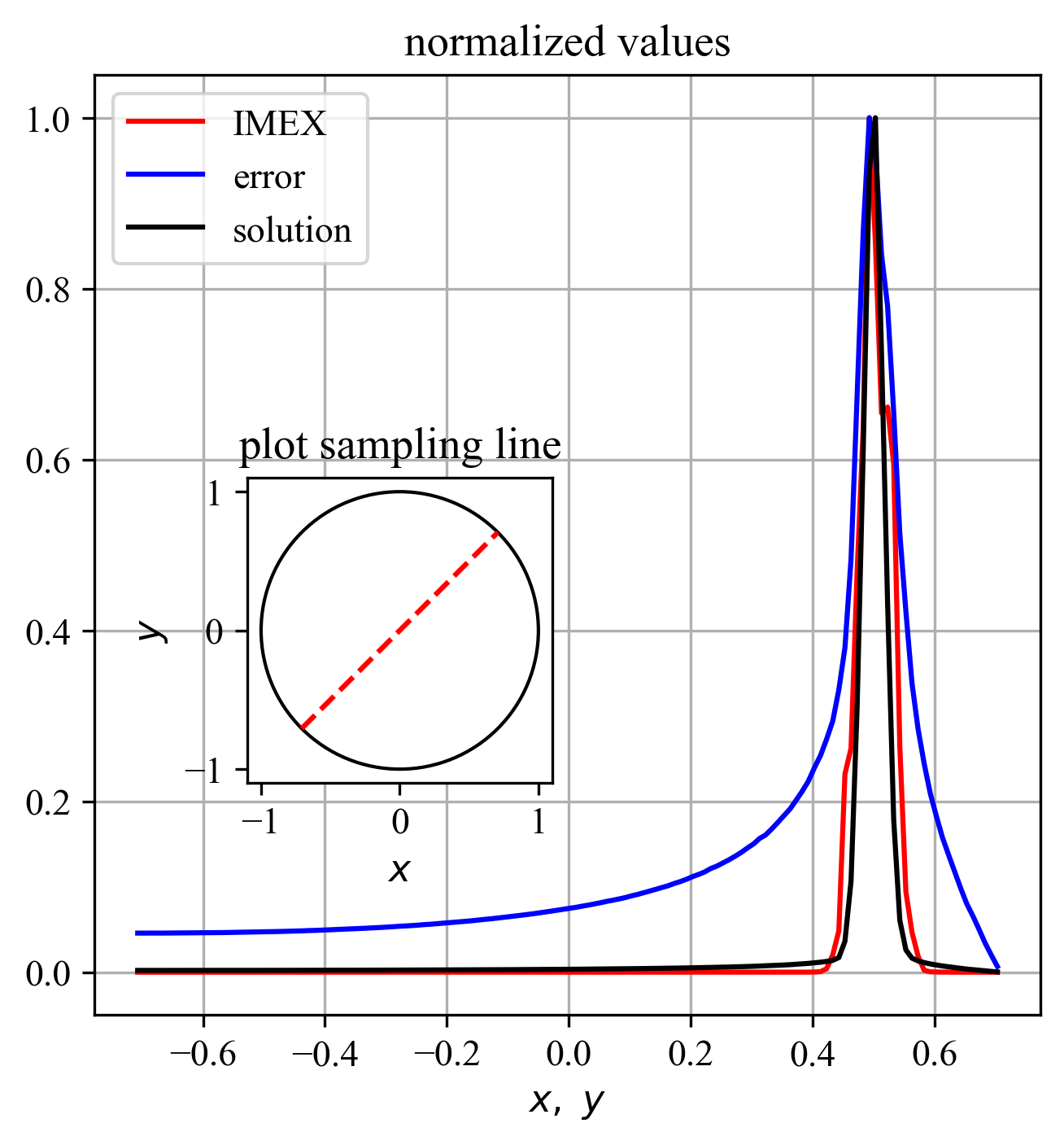}
	\caption{Implicit solution (solution), $\epsilon_{an}$ (error), and $\epsilon_{IMEX}$ (IMEX) normalized to their respective maximal values, plotted along the line $y = x; \ x,y \in \Omega$, $\L^{(hi)}$ are computed with 6\textsuperscript{th} degree monomials.}\label{fig: lines6}
\end{figure}
Additionally, the same case is solved with 6\textsuperscript{th} degree monomials used to produce $\L^{(hi)}$ for IMEX, with results plotted in Fig.~\ref{fig: lines6}.

Comparing Figs.~\ref{fig: solution}, and \ref{fig: imex} it is noticeable that the solution's error is the biggest around the source at point $\x _s =(0.5, 0.5)$.
The IMEX error indicator also predicts the biggest error to be around the same point, as can be seen in Fig.~\ref{fig: imex}.
This is further supported by the graph in Fig.~\ref{fig: lines}.
Although the IMEX error indicator does not follow the actual error, it successfully identifies the area of the biggest error.
Increasing the monomial degree to compute $\L^{(hi)}$ does not noticeably impact IMEX's performance, as can be seen by comparing Fig.~\ref{fig: lines}, and \ref{fig: lines6}.
However, increasing the monomial degree results in a significant compute performance hit in this particular case (total computation time increased to 4 s, compared to previous 2 s).

\section{Conclusions}
A synthetic example of the Poisson equation was solved and the IMEX error indicator was tested on it.
The error indicator correctly indicated the area of increased error, which also coincided with the source in the Poisson equation.
Results were produced by increasing the monomial degree of the explicit approximations by 2 compared to the implicit counterparts. Further increasing the monomial degree did not prove beneficial in this specific example.

We show that the proposed error indicator successfully identifies the areas with high error of the numerical solution. In the continuation, these findings could be used to adaptively refine the critical areas and improve the precision of the numerical solution.

\bibliographystyle{IEEEtran}
\bibliography{refs}

\begin{thebibliography}{10}
\providecommand{\url}[1]{#1}
\csname url@samestyle\endcsname
\providecommand{\newblock}{\relax}
\providecommand{\bibinfo}[2]{#2}
\providecommand{\BIBentrySTDinterwordspacing}{\spaceskip=0pt\relax}
\providecommand{\BIBentryALTinterwordstretchfactor}{4}
\providecommand{\BIBentryALTinterwordspacing}{\spaceskip=\fontdimen2\font plus
\BIBentryALTinterwordstretchfactor\fontdimen3\font minus
  \fontdimen4\font\relax}
\providecommand{\BIBforeignlanguage}[2]{{%
\expandafter\ifx\csname l@#1\endcsname\relax
\typeout{** WARNING: IEEEtran.bst: No hyphenation pattern has been}%
\typeout{** loaded for the language `#1'. Using the pattern for}%
\typeout{** the default language instead.}%
\else
\language=\csname l@#1\endcsname
\fi
#2}}
\providecommand{\BIBdecl}{\relax}
\BIBdecl

\bibitem{tolstykh2000using}
A.~I. Tolstykh, ``On using rbf-based differencing formulas for unstructured and
  mixed structured-unstructured grid calculations,'' in \emph{Proceedings of
  the 16th IMACS world congress}, vol. 228.\hskip 1em plus 0.5em minus
  0.4em\relax Lausanne, 2000, pp. 4606--4624.

\bibitem{tolstykh2003using}
A.~Tolstykh and D.~Shirobokov, ``On using radial basis functions in a “finite
  difference mode” with applications to elasticity problems,''
  \emph{Computational Mechanics}, vol.~33, no.~1, pp. 68--79, 2003.

\bibitem{belytschko1996meshless}
T.~Belytschko, Y.~Krongauz, D.~Organ, M.~Fleming, and P.~Krysl, ``Meshless
  methods: an overview and recent developments,'' \emph{Computer methods in
  applied mechanics and engineering}, vol. 139, no. 1-4, pp. 3--47, 1996.

\bibitem{bayona2017role}
V.~Bayona, N.~Flyer, B.~Fornberg, and G.~A. Barnett, ``On the role of
  polynomials in rbf-fd approximations: Ii. numerical solution of elliptic
  pdes,'' \emph{Journal of Computational Physics}, vol. 332, pp. 257--273,
  2017.

\bibitem{flyer2016}
\BIBentryALTinterwordspacing
N.~Flyer, B.~Fornberg, V.~Bayona, and G.~A. Barnett, ``On the role of
  polynomials in rbf-fd approximations: I. interpolation and accuracy,''
  \emph{Journal of Computational Physics}, vol. 321, pp. 21--38, 2016.
  [Online]. Available:
  \url{https://www.sciencedirect.com/science/article/pii/S0021999116301632}
\BIBentrySTDinterwordspacing

\bibitem{9597066}
J.~Slak, ``Partition-of-unity based error indicator for local collocation
  meshless methods,'' in \emph{2021 44th International Convention on
  Information, Communication and Electronic Technology (MIPRO)}, 2021, pp.
  254--258.

\bibitem{carstensen2005explicit}
C.~Carstensen, R.~Lazarov, and S.~Tomov, ``Explicit and averaging a posteriori
  error estimates for adaptive finite volume methods,'' \emph{SIAM Journal on
  Numerical Analysis}, vol.~42, no.~6, pp. 2496--2521, 2005.

\bibitem{SEGETH20101589}
\BIBentryALTinterwordspacing
K.~Segeth, ``A review of some a posteriori error estimates for adaptive finite
  element methods,'' \emph{Mathematics and Computers in Simulation}, vol.~80,
  no.~8, pp. 1589--1600, 2010, eSCO 2008 Conference. [Online]. Available:
  \url{https://www.sciencedirect.com/science/article/pii/S0378475408004230}
\BIBentrySTDinterwordspacing

\bibitem{slak2021}
J.~Slak and G.~Kosec, ``{Medusa: A C++ Library for solving PDEs using Strong
  Form Mesh-Free methods},'' \emph{{ACM Transactions on Mathematical
  Software}}, 2021.

\bibitem{rot}
\BIBentryALTinterwordspacing
M.~Rot and G.~Kosec, ``Refined rbf-fd analysis of non-newtonian natural
  convection,'' 2022. [Online]. Available:
  \url{https://arxiv.org/abs/2202.08095}
\BIBentrySTDinterwordspacing

\bibitem{filip}
\BIBentryALTinterwordspacing
F.~Strniša, M.~Jančič, and G.~Kosec, ``A meshless solution of a small-strain
  plasticity problem,'' 2022. [Online]. Available:
  \url{https://arxiv.org/abs/2203.08462}
\BIBentrySTDinterwordspacing

\bibitem{mitchell2013}
\BIBentryALTinterwordspacing
W.~F. Mitchell, ``A collection of 2d elliptic problems for testing adaptive
  grid refinement algorithms,'' \emph{Appl. Math. Comput.}, vol. 220, p.
  350–364, sep 2013. [Online]. Available:
  \url{https://doi.org/10.1016/j.amc.2013.05.068}
\BIBentrySTDinterwordspacing

\bibitem{slak2019}
J.~Slak and G.~Kosec, ``On generation of node distributions for meshless pde
  discretizations,'' \emph{SIAM Journal on Scientific Computing}, vol.~41,
  no.~5, pp. A3202--A3229, 2019.

\bibitem{eigenweb}
G.~Guennebaud, B.~Jacob \emph{et~al.}, ``Eigen v3,''
  http://eigen.tuxfamily.org, 2010.

\bibitem{larry_bradley_2020_4044744}
\BIBentryALTinterwordspacing
L.~Bradley, B.~Sip{\H o}cz, T.~Robitaille, E.~Tollerud, Z.~Vin{\'{\i}}cius,
  C.~Deil, K.~Barbary, T.~J. Wilson, I.~Busko, H.~M. G{\"u}nther, M.~Cara,
  S.~Conseil, A.~Bostroem, M.~Droettboom, E.~M. Bray, L.~A. Bratholm, P.~L.
  Lim, G.~Barentsen, M.~Craig, S.~Pascual, G.~Perren, J.~Greco, A.~Donath,
  M.~de~Val-Borro, W.~Kerzendorf, Y.~P. Bach, B.~A. Weaver, F.~D'Eugenio,
  H.~Souchereau, and L.~Ferreira, ``astropy/photutils: 1.0.0,'' Sep. 2020.
  [Online]. Available: \url{https://doi.org/10.5281/zenodo.4044744}
\BIBentrySTDinterwordspacing

\end{thebibliography}

\end{document}